\documentclass{article}

\usepackage{amsmath,amsfonts,amssymb,amsthm,fullpage}

\newtheorem{theorem}{Theorem}[section]

\theoremstyle{remark}

\begin{document}

\title{Three--dimensional Gaussian fluctuations of spectra of overlapping stochastic Wishart matrices}

\author{Jeffrey Kuan and Zhengye Zhou}

\date{}

\maketitle

\abstract{In \cite{DP18}, the authors consider eigenvalues of overlapping Wishart matrices and prove that its fluctuations asymptotically convergence to the Gaussian free field. In this brief note, their result is extended to show that when the matrix entries undergo stochastic evolution, the fluctuations asymptotically converge to a three-dimensional Gaussian field, which has an explicit contour integral formula. This is analogous to the result of \cite{Bo3D1} for stochastic Wigner matrices.}

\section{Introduction and Main Result}
\subsection{Motivation}
The two--dimensional Gaussian free field appears as the asymptotic fluctuations of the height functions of several two--dimensional models, such as eigenvalues of sub--matrices of Wigner random matrices and interacting particle systems on a two--dimensional lattice. For example, in \cite{BFCMP}, the authors prove that an interacting particle system on $\mathbb{Z} \times \mathbb{Z}_+$  has asymptotic fluctuations described by the Gaussian free field on the upper half--plane $\mathbb{H}$. This interacting particle system can be constructed through the representation theory of the infinite--dimensional unitary group $U(\infty)$ \cite{BK08}. Similarly, in \cite{Bo3D1} it is proven that for Wigner matrices whose entries undergo stochastic evolution, the fluctuations of the eigenvalues of sub--matrices are described by a three--dimensional Gaussian field which generalizes the Gaussian free field on $\mathbb{H}$. This same three--dimensional field occurs in a non--commutative random walk \cite{K1}.

In \cite{BK10}, the representation theory of the infinite--dimensional orthogonal group $O(\infty)$ is used to construct an interacting particle system, this time on the quadrant $\mathbb{Z}_{ \geq 0} \times \mathbb{Z}_+$. The vertical boundary $\{0\} \times \mathbb{Z}_+$ acts as a reflecting barrier for the particles. In \cite{KuanGFF}, it is proven that the fixed--time asymptotic fluctuations of the interacting particle system are the Gaussian free field on the upper half--plane with the unit disk removed, $\mathbb{H}-\mathbb{D}$. Upcoming work by the second author generalizes this result to multiple time parameters, wherein a three--dimensional Gaussian field appears in the asymptotic limit. Thus, it is natural to ask if the same three--dimensional Gaussian field occurs in the context of random matrices. An appropriate analog is Wishart matrices, because the eigenvalues are always non--negative, corresponding to the interacting particles being restricted to be to the right of the reflecting barrier. In \cite{DP18}, it is shown that the eigenvalues of overlapping Wishart matrices have fluctuations that converge to the Gaussian free field. In the present work, their latter result is generalized to the case when the entries of the Wishart matrices undergo stochastic evolution. 

\subsection{Main Result}
Consider an infinite array $\{Z_{ij}(t):i,j\geq 1\}$ of stochastic processes (not necessarily Markov). The $Z_{ij}(t)$ are allowed to be real, complex, or quaternionic, corresponding to the usual $\beta=1,2,4$ values. In all cases, assume that all real and imaginary components of $Z_{ij}(t)$ are mutually independent. Assume the following moment conditions for all fixed times $t\geq 0$:
$$
\mathbb{E}Z_{ij}(t)=0, \quad \mathbb{E}\vert Z_{ij}(t)\vert^2 = 1, \quad \mathbb{E}\vert Z_{ij}(t)\vert^4= 1 + \frac{2}{\beta}, \quad \sup_{i,j\geq 1} \mathbb{E}\vert Z_{ij}(t)\vert^{4+\epsilon}<\infty.
$$
The first, second, and fourth moments agree with the standard real, complex and quaternion Gaussian random variables for all $t\geq 0$. Further suppose that for all $s,t\geq 0$
$$
\mathbb{E} \left|  Z_{ij}(s) \right|^2  \left|Z_{ij}(t)\right|^2 = c_{2}(s,t), \quad \quad \mathbb{E} \vert Z_{ij}(s)Z_{ij}(t)] \vert = c_1(s,t).
$$

From this infinite array, consider rectangular submatrices $B_i(s)$ for $1\leq i \leq k$ with sizes $m_iL \times n_iL$ that overlap at $m_{ij}L$ rows and $n_{ij}L$ columns. Define the Wishart matrices $W_i(s)=B_i(s)^*B_i(s)/L$. The statistics of interest are
$$
X_i^{(p)}(s) := \mathrm{Tr}(W_i^p(s))) - \mathbb{E} \mathrm{Tr}(W_i^p(s)) 
$$
for positive integers $p$. The result of this paper is that as $L\rightarrow \infty$, the random variables $X_i^{(p)}$ converge to a Gaussian field with an explicit covariance.  

\begin{theorem}
Suppose that $m_i,n_j,m_{ij},n_{ij}$ all depend on $L$ in such a way that
$$
\frac{m_i}{L}\rightarrow \mu_i, \quad \quad  \frac{n_j}{L} \rightarrow \nu_j, \quad \quad \frac{m_{ij}n_{ij}L^2}{m_in_im_jn_j}\rightarrow \theta_{ij}
$$
as $L\rightarrow \infty$. Then $(X_1^{(p_1)}(s_1), \ldots, X_k^{(p_k}(s_k))$ converges in distribution to a centered Gaussian vector $(\xi_1,\ldots,\xi_k)$ with covariance
\begin{multline*}
\mathbb{E}[\xi_i\xi_j] = \frac{4c_1(s_i,s_j)p_ip_j}{\beta (i \pi)^2 }\oint\oint  (\mu_i + \nu_i + 2\Re \zeta_i)^{p_i-1} (\mu_j+ \nu_j + 2\Re \zeta_j)^{p_i-1} \log \left| \frac{\theta_{ij}^{-1} - {\zeta_i}\zeta_j}{\theta_{ij}^{-1} - {\zeta_i}\bar{\zeta_j}} \right| \frac{\Im \zeta_i}{\zeta_i}  \frac{\Im \zeta_j}{\zeta_j} d \zeta_i d\zeta_j\\
+ \frac{ 4(c_2(s_i,s_j)-1-2\beta^{-1}c_1(s_i,s_j))  \theta_{ij}p_ip_j   }{ ( \pi i)^2 } \oint\oint (\mu_i + \nu_i + 2\Re \zeta_i)^{p_i-1} (\mu_j+ \nu_j + 2\Re \zeta_j)^{p_i-1}  \frac{\left( \Im \zeta_i\right)^2}{\zeta_i}  \frac{\left( \Im \zeta_j\right)^2}{\zeta_j} d \zeta_i d\zeta_j
\end{multline*}
where the $\zeta_i,\zeta_j$ contours are semicircles over the upper half plane $\mathbb{H}$ with radii $\sqrt{\mu_i\nu_i},\sqrt{\mu_j\nu_j}., $ respectively.
\end{theorem}

\textbf{Remark}. When $s_i=s_j$, this reduces to Lemma 2.2 of \cite{DP18}. Also, by the Cauchy--Schwarz inequality, $c_2(s,t) \leq \sqrt{ \mathbb{E}| Z_{ij}(s)|^4 |Z_{ij}(t)|^4} = 1+2/\beta$, so the constant in the second line is always non--positive.

\textbf{Acknowledgments}. The authors would like to thank Alexei Borodin, Ioana Dumitriu and Elliot Paquette for helpful discussions. The simulations done in section 3 were aided by Texas A\&M University High Performance Research Computing.

\section{Proof}
The proof is in essence identical to \cite{DP18}. As explained in \cite{DP18}, the centered trace $X^{(p)}_i(t)$ can be written as a sum over closed walks of length $2p$ on the complete bipartite graph $K(\mathbb{N},\mathbb{N}). $ More explicitly, represent a walk as a function $w:\{1,\ldots,2p\} \rightarrow \mathbb{N}$ and define
$$
Z_w(t) = \prod_{i=1}^p Z_{w(2p-1),w(2p)}(t) \overline{Z_{w(2p-1),w(2p)}(t) } - \mathbb{E}\left[   \prod_{i=1}^p Z_{w(2p-1),w(2p)}(t) \overline{Z_{w(2p-1),w(2p)}(t) }  \right].
$$
With this definition, 
$$
X_i^{(p)}(t) = \sum_{w \in \mathcal{S}_p} Z_w(t),
$$
where $\mathcal{S}_p$ is the set of all closed walks of length $2p$ on the complete bipartite graph $K(\mathbb{N},\mathbb{N})$.

Thus, when computing the covariance of $X_i^{(k)}$ and $X_j^{(l)}$, the sum is over pairs of walks. As explained on page 8 of \cite{DP18}, the only pairs of walks that contribute in the limit are those that cover exactly $k+l$ vertices and traverse a common edge. There are two types of pairs of walks that can arise this way. The first type is that both walks are depth--first search walks of two trees that cover a common edge. For example, the first walk can be $1\rightarrow 2 \rightarrow 1 \rightarrow 3 \rightarrow 1$ while the second walk is $1\rightarrow 2 \rightarrow 1 \rightarrow 4 \rightarrow 1$. The contributions of these walks is then
$$
 \frac{ 4(\mathbb{E}[Z_{11}^2(s)Z_{11}^2(t)]-1) k\cdot l \cdot \theta_{ij}}{( \pi i)^2} \oint\oint\ (\mu_1 + \nu_1 + 2\Re\zeta_1)^{k-1} (\mu_2 + \nu_2 + 2\Re\zeta_2)^{l-1}\frac{(\Im \zeta_1)^2}{\zeta_1}\frac{(\Im \zeta_2)^2}{\zeta_2} d\zeta_1 d\zeta_2.
$$
(compare to equation (8)) of \cite{DP18}. The second type of pairs of walks are when both walks cover a unicyclic graph, traversing the cycle once and making excursions along trees attached to this cycle. The cycle must be common to both walks so that it is traversed twice by the union of the walks. This type of walk contributes
\begin{multline*}
 -\frac{8\mathbb{E}\vert Z_{ij}(s)Z_{ij}(t)\vert k \cdot  l \cdot  \theta_{ij}}{\beta( \pi i)^2}  \oint\oint\ (\mu_1 + \nu_1 + 2\Re\zeta_1)^{k-1} (\mu_2 + \nu_2 + 2\Re\zeta_2)^{l-1} \frac{(\Im \zeta_1)^2}{\zeta_1}\frac{(\Im \zeta_2)^2}{\zeta_2} d\zeta_1 d\zeta_2 \\
 \frac{4\mathbb{E}[Z_{ij}(s)Z_{ij}(t)]kl}{\beta(\pi i)^2}\oint\oint  (\mu_1 + \nu_1 + 2\Re\zeta_1)^{k-1} (\mu_2 + \nu_2 + 2\Re\zeta_2)^{l-1}  \log \left| \frac{\theta_{ij}^{-1} - {\zeta_1}\zeta_2}{\theta_{ij}^{-1} - {\zeta_1}\bar{\zeta_2}}\right| \frac{\Im \zeta_1}{\zeta_1}\frac{\Im \zeta_2}{\zeta_2} d\zeta_1 d\zeta_2 
\end{multline*}
Combining the two sums yields the expression in the statement of the Theorem.

Finally, the asymptotic fluctuations are asymptotically Gaussian because of the same arguments as in section 3 of \cite{DP18}.

\section{Examples}
Integrating by parts twice yields, for the main term,
$$
\frac{c_1(s_i,s_j)\theta_{ij}^{-1}}{2\beta \pi^2} \oint\oint \frac{ (\mu_i + \nu_i + 2\Re(\zeta_i))^{p_i}  (\mu_j + \nu_j+ 2\Re(\zeta_j))^{p_j} }{ (\theta_{ij}^{-1} - \zeta_i\zeta_j)^2}d \zeta_i d\zeta_j
$$
where the $\zeta_i,\zeta_j$ contours are now full circles in $\mathbb{C}$ with radii $\sqrt{\mu_i\nu_i},\sqrt{\mu_j\nu_j}, $ respectively. The error term can be evaluated directly with residues without integrating by parts.  For example, if $p_i=p_j=1$ then, using the Taylor series $(\theta_{ij}^{-1}-\zeta_i\zeta_j)^{-2} = \theta_{ij}^2(1+2\theta_{ij}^{} \zeta_i^{}\zeta_j^{}+\cdots),$ the covariance equals
$$
2c_1(s_i,s_j)\beta^{-1}\theta_{ij}\mu_i\nu_i\mu_j\nu_j + \theta_{ij}\mu_i\mu_j\nu_i\nu_j(c_2(s_i,s_j)-1-2\beta^{-1}c_1(s_i,s_j)).
$$
This can be checked with Python: for example, if $Z_{11}(t)$ is a real Ornstein--Uhlenbeck process with $\mathbb{E}[Z_{11}(s_i)Z_{11}(s_j)]=1/2$, then $\mathbb{E}[Z_{11}(s_i)^2Z_{11}(s_j)^2]=1+2\cdot (1/2)^2 = 3/2$, and if $\mu_i=4,\nu_i=2,\mu_j=\nu_j=1=m_{ij}=n_{ij}$, then $\theta_{ij}=1/8$ and the prediction for the covariance is $2\cdot 1/2\cdot 1/8 \cdot 8 + 1/8 \cdot 8  \cdot (3/2-1-2\cdot 1/2)=1/2$, and a simulation of 16,800 samples of Wishart matrices at $L=400$ yields $\approx 0.51310$ for the covariance. Higher moments are difficult to simulate, because the Gaussian free field has logarithmic fluctuations, requiring essentially exponential computation time. However, with the aid of Texas A\&M High Performance Research Computing, 48,000 samples at size $L=400$ predicts a covariance of $\approx181$ when $p_i=p_j=2$ and $(\mu_i,\nu_i,\mu_j,\nu_j, m_{ij},n_{ij})=(5,2,3,1,1,1) ,$ and $Z_{11}(t)$ as above. The formula yields
\begin{multline*}
2c_1(s_i,s_j)\theta_{ij} \mu_i\nu_i\mu_j\nu_j\left(  4(\mu_i+\nu_i)(\mu_j+\nu_j) 
+ 2\theta_{ij} \mu_i\nu_i\mu_j\nu_j \right) \\
+ (c_2(s_i,s_j)-1-2\beta^{-1}c_1(s_i,s_j))\theta_{ij}\mu_i\nu_i\mu_j\nu_j\cdot 4(\mu_i+\nu_i)(\mu_j+\nu_j) ,
\end{multline*}
which numerically evaluates to $354-168=186$. Both simulations are within $3\%$ of the correct answer.

The Python code is available from the first author upon request.

\bibliographystyle{alpha}
\bibliography{Wishart}

\end{document}